\documentclass{ACmod}

\DeclareFontFamily{U}{rsf}{}
\DeclareFontShape{U}{rsf}{m}{n}{
  <5> <6> rsfs5 <7> <8> <9> rsfs7 <10-> rsfs10}{}
\DeclareMathAlphabet{\mathscr}{U}{rsf}{m}{n}

\DeclareMathAlphabet{\mathgth}{U}{euf}{m}{n}

\DeclareFontFamily{U}{cyr}{}
\DeclareFontShape{U}{cyr}{m}{n}{
  <5> wncyr5 <6> wncyr6 <7> wncyr7 <8> wncyr8 <9> wncyr9 <10-> wncyr10}{}
\DeclareMathAlphabet{\mathcyr}{U}{cyr}{m}{n}

\input cyracc.def

\DeclareSymbolFont{bbold}{U}{bbold}{m}{n}
\DeclareSymbolFontAlphabet{\mathbbold}{bbold}

\makeatletter
\def\operator@font{\sf}
\makeatother

\renewcommand{\phi}{\varphi}

\newcommand{\F}{\mathcal{F}}

\usepackage{fancyhdr}
\usepackage{amssymb}
\usepackage{xypic}

\usepackage{tikz-cd}
\usepackage{mathtools}
\usepackage{amsfonts}

\usepackage{tikz}
\usepackage[tocflat]{tocstyle}
\usetikzlibrary{positioning}
\usepackage{amscd}

\usepackage{amsthm}
\usepackage{amsmath}

\everymath{\displaystyle}

\begin{document}

\title{A note on a question of Markman}

\author[Shengyuan Huang]{%
Shengyuan Huang}

\address{Shengyuan Huang, Mathematics Department,
University of Wisconsin--Madison, 480 Lincoln Drive, Madison, WI
53706--1388, USA.\\
Email: \href{mailto:shuang279@wisc.edu}{shuang279@wisc.edu}
}

\begin{abstract}
  {\sc Abstract:} Let $\F$ be a vector bundle on a complex projective algebraic variety $X$.  If $\F$ deforms along a first order deformation of $X$, its Mukai vector remains of Hodge type along this deformation.  We prove an analogous statement for all polyvector fields, not only for those in $H^1(X,T_X)$ corresponding to deformations of the complex structure.  This answers a question of Markman.  We also explore a Lie theoretic analogue of the statement above.\\


\end{abstract}

\maketitle

\setcounter{tocdepth}{1}
\tableofcontents

\section{Introduction}

\paragraph\label{1.1}
Let $X$ be a smooth complex algebraic variety.  Consider the first order deformation $\tilde{X}$ of $X$ associated to a class $\tilde{\alpha}\in H^1(X,T_X)$.

In general, a vector bundle $\F$ may not deform to a bundle $\tilde{\F}$ on $\tilde{X}$.  The obstruction $\alpha_\F\in\mathrm{Ext}^2(\F,\F)$ to the existence of a vector bundle $\tilde{\F}$ on $\tilde{X}$ such that $\tilde{\F}|_X\cong\F$ was described in \cite{B,T} as the contraction
$$\alpha_\F=\tilde{\alpha}\lrcorner\, at_\F\in\mathrm{Ext}^2(\F,\F).$$
Here $at_\F\in\mathrm{Ext}^1(\F,\F\otimes\Omega_X)$ is the Atiyah class of $\F$.

Moreover, if $\F$ does deform, then its Chern classes and hence its Mukai vector stay of Hodge type on the deformed space $\tilde{X}$.  This implies that the class
$$\tilde{\alpha}\lrcorner\, v(\F)\in\mathrm{H\Omega}_*(X)\overset{\mathrm{def}}{=}\bigoplus_{q-p=*}H^p(X,\wedge^q\Omega_X)$$
vanishes, where $v(\F)$ is the Mukai vector of $\F$.

Thus, in the simple case where $\tilde{\alpha}\in H^1(X,T_X)$ we conclude that if $\tilde{\alpha}\lrcorner\,at_\F$ is zero, then $\tilde{\alpha}\lrcorner\,v(\F)$ is zero.

\paragraph In an email correspondence, Eyal Markman asked if the above statement can be generalized to the case where $\tilde{\alpha}$ is an arbitrary polyvector field in $\mathrm{HT}^*(X)=\bigoplus_{p+q=*}H^p(X,\wedge^q T_X)$.  According to Markman, this question is central to his study of the deformations of hyperk\"{a}hler manifolds.  We provide an answer to this question in this paper.

\paragraph First, there are two ways to generalize the class that appears in~(\ref{1.1}) above.  For a class $\tilde{\alpha}\in\mathrm{HT}^*(X)$ we can define classes in $\mathrm{Ext}^*(\F,\F)$ in the following two ways.

The first is defined by using the HKR isomorphism
$$I^{HKR}:\mathrm{HT}^*(X)\rightarrow\mathrm{HH}^*(X),$$
where the latter is the Hochschild cohomology of $X$.  Since $\mathrm{HH}^*(X)$ can be interpreted as natural transformations of the identity functor at the dg level, this yields a natural map $\mathrm{HH}^*(X)\rightarrow\mathrm{Ext}^*(\F,\F)$.

The second construction was defined by Toda \cite{T}.  Consider the exponential Atiyah class $$\exp(at_\F)=1+at_\F+\cdots+\frac{(at_\F)^k}{k!}+\cdots,$$
where $(at_\F)^k\in\mathrm{Ext}^k(\F,\F\otimes\wedge^k\Omega_X)$.  Let $\tilde{\alpha}^{p,k}\in H^p(X,\wedge^kT_X)$ be the homogenous degree $(p,k)$ part of $\tilde{\alpha}$.  We can contract $\tilde{\alpha}^{p,k}$ with $\frac{(at_\F)^k}{k!}$ to get an element in $\mathrm{Ext}^{p+k}(\F,\F)$.  Taking the sum over all $(p,k)$, we get the desired class which will be denoted by $\tilde{\alpha\lrcorner\exp(at_\F)}\in\mathrm{Ext}^*(\F,\F)$.  When $\tilde{\alpha}$ is a class in $H^1(X,T_X)$, we recover the previous contraction $\tilde{\alpha}\lrcorner\,at_\F$.

Our first result is below.
\paragraph{\bf Theorem A.}
{\em
The two classes defined above are the same.  In other words the diagram
$$\xymatrix{
\mathrm{HH}^*(X)\ar[r] & \mathrm{Ext}^*(\F,\F)\\
\mathrm{HT}^*(X)\ar[u]^{I^{HKR}}\ar[ur]_{(-)\lrcorner\exp(at_\F)}
}
$$
is commutative.}

There is an analogous result for Hopf algebras.  See Theorem 2.7 in \cite{CI} and see \cite{PW} for more details.

\paragraph The space $\mathrm{H\Omega}_*(X)$ is naturally a module over $\mathrm{HT}^*(X)$, mimicking the module structure of Hochschild homology over cohomology.  For an object $\F$ in the derived category of $X$, its Mukai vector $v(\F)$ lies in $\mathrm{H\Omega}_*(X)$.  Thus we can act with the class $\tilde{\alpha}$ to obtain $\tilde{\alpha}\lrcorner\,v(\F)\in\mathrm{H\Omega}_*(X)$.

{\bf Theorem B.}
{\em  If $\tilde{\alpha}\lrcorner\exp(at_\F)=0$, then we have
$$D(\tilde{\alpha})\lrcorner\,v(\F)=0.$$
Here $D$ is the Duflo operator,
 $$D(\tilde{\alpha})=Td(X)^{\frac{1}{2}}\lrcorner\,\tilde{\alpha},$$
 where $Td(X)$ is the Todd class of $X$.
}

\begin{remark}
We are using the contraction symbol $\lrcorner$ in three different ways in this paper.
\begin{itemize}
\item A polyvector field $\tilde{\alpha}\in\mathrm{HT}^*(X)$ acts on a class $v\in\mathrm{H\Omega}_*(X)$.  This action is denoted by $\tilde{\alpha}\lrcorner\,v\in\mathrm{H\Omega}_*(X)$.
\item A class $v\in\mathrm{H\Omega}_*(X)$ acts on a polyvector field $\tilde{\alpha}\in\mathrm{HT}^*(X)$.  This action yields an element $v\lrcorner\,\tilde{\alpha}\in\mathrm{HT}^*(X)$.  We only use the second contraction in the Duflo operator $D(\tilde{\alpha})=Td(X)^{\frac{1}{2}}\lrcorner\,\tilde{\alpha}$ in this paper.  Note that $D$ is an automorphism of $\mathrm{HT}^*(X)$.  The inverse operator is $D^{-1}(\tilde{\alpha})=Td(X)^{-\frac{1}{2}}\lrcorner\,\tilde{\alpha}$.
\item The third contraction map is $\beta\lrcorner\exp(at_\F)\in\mathrm{Ext}^*(\F,\F)$ for $\beta\in\mathrm{HT}^*(X)$.  An element $\beta$ in $H^p(X,\wedge^k T_X)$ can only contract with the term $\frac{(at_\F)^k}{k!}$ in the Taylor expansion of $\exp(at_\F)$.  It is easy to distinguish this map from the previous two maps.
\end{itemize}
\end{remark}
\paragraph The inspiration for Theorem A comes from a similar statement in Lie theory.  Let $\mathfrak{g}$ be a finite dimensional Lie algebra and $V$ be a finite dimensional representation.  One can draw the diagram
$$
\xymatrix{
(U\mathfrak{\mathfrak{g}})^{\mathfrak{g}}\ar[r] & \mathrm{Hom}(V,V)\\
(S\mathfrak{g})^{\mathfrak{g}}\ar[u]\ar[ur]
}
$$
which is similar to the one in Theorem A.  We will provide more details and prove that the diagram above is commutative in section 2.

\paragraph Note that our statement in Theorem B appears to be different from the original one, which did not have the Duflo operator $D$.  We will prove that the original statement follows easily from ours.

\paragraph {\bf Plan of the paper.} Section 2 contains the proof of Theorem A and of its Lie theoretic analogue.

Section 3 is devoted to the proof of Theorem B.  It is a consequence of Theorem A.  At the end we prove that we can recover the result in~(\ref{1.1}) from Theorem B.

\paragraph {\bf Acknowledgments.} I would like to thank Andrei C\u{a}ld\u{a}raru for discussing details with me during our weekly meeting.  I have benefited from stimulating email correspondence with Dror Bar-Natan.  I am also grateful to Eyal Markman for many valuable comments.  The author is partially supported by the National Science Foundation under Grant No. DMS-1811925.

\section{The proof of Theorem A}
We prove Theorem A in this section.  The diagram in Theorem A has a Lie theoretic background.  We draw the corresponding diagram for Lie algebras and we explain the similarity between the Lie algebra diagram and the diagram in Theorem A.  We provide a proof for the commutativity of the Lie algebra diagram and explain that the proof can be generalized to the diagram in Theorem A.

\paragraph {\bf A similar diagram for Lie algebras.} \label{Lie diagram}Let $\mathfrak{g}$ be a finite dimensional Lie algebra over a field of characteristic zero and let $V$ be a finite dimensional representation of $\mathfrak{g}$.  There is a diagram
$$
\xymatrix{
(U\mathfrak{g})^{\mathfrak{g}}\ar[r] & \mathrm{Hom}(V,V)\\
(S\mathfrak{g})^{\mathfrak{g}}.\ar[u]^{\mbox{PBW}}\ar[ur]
}
$$
The PBW map from the symmetric algebra $S\mathfrak{g}$ to the universal enveloping algebra $U\mathfrak{g}$ is defined on the degree $n$-th component of $S\mathfrak{g}$ as follows
$$x_1\cdots x_n\rightarrow\frac{1}{n!}\sum_{\sigma\in S_n}x_{\sigma(1)}\cdots x_{\sigma(n)}.$$
Here $S_{n}$ is the symmetric group on a finite set of n symbols.  The universal enveloping algebra $U\mathfrak{g}$ acts naturally on $V$.  This natural action defines the map $(\mathrm{U}\mathfrak{g})^{\mathfrak{g}}\rightarrow \mathrm{Hom}(V,V)$ on the top of the diagram above.  The map $(\mathrm{S}\mathfrak{g})^{\mathfrak{g}}\rightarrow \mathrm{Hom}(V,V)$ is defined as follows.  We can rewrite the representation map $\mathfrak{g}\otimes V\rightarrow V$ as a map $\Lambda:V\rightarrow V\otimes\mathfrak{g}^*$.  Take the exponent
$$\exp(\Lambda)=id_V+\Lambda+\cdots+\frac{\Lambda^k}{k!}+\cdots$$
of the map $\Lambda$.  Then we can contract $\exp(\Lambda)$ with $S\mathfrak{g}$.

In algebraic geometry, Kapranov and Kontsevich \cite{K} observed that the shifted tangent bundle $T_X[-1]$ has a Lie algebra structure in the derived category of $X$. Roberts and Willerton \cite{RW} proved that the category of representations of $T_X[-1]$ is the derived category of $X$ and the universal enveloping algebra of $T_X[-1]$ is the  Hochschild cochain complex $\mathcal{RH}om(\Delta_*\mathcal{O}_X,\Delta_*\mathcal{O}_X)$, where $\Delta$ is the diagonal embedding $\Delta:X\hookrightarrow X\times X$.  The functor $(-)^{\mathfrak{g}}$ is the 0-th Lie algebra cohomology which is similar to $H^*(X,-)$.  Setting $\mathfrak{g}$ to be equal to $T_X[-1]$ in the Lie algebra diagram, we get the diagram in Theorem A for a smooth complex variety $X$.  The Hochschild cohomology $\mathrm{HH}^*(X)$ plays the role of $(U\mathfrak{g})^{\mathfrak{g}}$, $\mathrm{HT}^*(X)$ plays the role of $(S\mathfrak{g})^{\mathfrak{g}}$, and the HKR map is precisely the PBW map.

\begin{Proof}[Proof of the commutativity for the Lie algebra diagram.]
We can prove that the diagram in~(\ref{Lie diagram}) is commutative even before taking $\mathfrak{g}$-invariants, i.e., the diagram
$$
\xymatrix{
U\mathfrak{g}\ar[r] & \mathrm{Hom}(V,V)\\
S\mathfrak{g}\ar[u]^{\mbox{PBW}}\ar[ur]
}
$$
is commutative.  The map PBW factors through the tensor algebra $T\mathfrak{g}$
$$\xymatrix{\mbox{PBW}:S\mathfrak{g}\ar[r]^{~~~~~~\psi} & T\mathfrak{g}\ar[r] &  U\mathfrak{g},}$$
so we can replace $U\mathfrak{g}$ at the top left corner of the diagram by $T\mathfrak{g}$.  It is easy to check that the map $S\mathfrak{g}\rightarrow\mathrm{Hom}(V,V)$ is equal to the following map
$$\xymatrix{ S\mathfrak{g}\ar[r]^{\psi} & T\mathfrak{g}\ar[r]^{\varphi~~~~~~~~} & \mathrm{Hom}(V,V),}$$
where the map $\varphi:T\mathfrak{g}\rightarrow\mathrm{Hom}(V,V)$ is defined as follows.  Rewrite the representation map $\mathfrak{g}\otimes V\rightarrow V$ as a map $\Lambda:V\rightarrow V\otimes\mathfrak{g}^*$.  Instead of taking the exponential of the map $\Lambda$, we compose the map $\Lambda$ with itself $k$ times.  We get a map $\Lambda^{\otimes k}:V\rightarrow V\otimes (\mathfrak{g}^*)^{\otimes k}$ in this way.  Contract $\Lambda^{\otimes k}$ with $\mathfrak{g}^{\otimes k}$ and get a map $\mathfrak{g}^{\otimes k}\rightarrow\mathrm{Hom}(V,V)$.  Adding the $k$-th components for all $k\in\mathbb{N}$, we obtain the desired map $\varphi: T\mathfrak{g}\rightarrow\mathrm{Hom}(V,V)$.

Now we have two maps $T\mathfrak{g}\rightarrow\mathrm{Hom}(V,V)$.  One of them is the map $\varphi$, and the other one is $\Theta: T\mathfrak{g}\rightarrow  U\mathfrak{g}\rightarrow \mathrm{Hom}(V,V)$.  We want to show that they agree. This follows from Lemma \ref{Lemma1} below by setting $W_1$ to be $V$ and $W_2$ to be $\mathfrak{g}^{\otimes k}$.

\end{Proof}

\begin{Lemma}\label{Lemma1}
Let $W_1$ and $W_2$ be finite dimensional vector spaces over a field $k$ and $f$ be a map $W_2\otimes W_1\rightarrow W_1$.  Rewrite the map as $g:W_1\rightarrow W_2^*\otimes W_1$ by the adjunction formula $\mathrm{Hom}(W_2\otimes_k W_1, W_1)=\mathrm{Hom}(W_1,W_2^*\otimes_kW_1)$.  Fix an element $x\in W_2$.  Then $f(x\otimes-)$ is a map from $W_1$ to $W_1$.  This map is precisely $g$ followed by the contraction with $x$.
\end{Lemma}
\begin{proof}
This is due to the adjunction property $$\mathrm{Hom}(W_2\otimes_k W_1, W_1)=\mathrm{Hom}(W_1,W_2^*\otimes_kW_1).$$
\end{proof}
\begin{Proof}[Proof of Theorem A.] The proof above reduces the commutativity of the Lie algebra diagram to a statement about tensor algebras.  The statement about tensor algebras remains valid in the case of derived categories.

One can define a map $\mathrm{Sym}(T_X[-1])\rightarrow T(T_X[-1])$ given by the formula
$$x_1\wedge\cdots \wedge x_n\rightarrow\frac{1}{n!}\sum_{\sigma\in S_n}(-1)^{sgn(\sigma)}x_{\sigma(1)}\cdots x_{\sigma(n)},$$
where $T(T_X[-1])$ is the tensor algebra on $T_X[-1]$.

The map above is a differential graded version of the map $\psi$ in~(\ref{Lie diagram}).  Let $X^{(1)}$ be the first order neighborhood of $X$ in $X\times X$.  There are embeddings $i:X\hookrightarrow X^{(1)}$ and $j:X^{(1)}\hookrightarrow X\times X$.  Arinkin and C\u{a}ld\u{a}raru \cite{AC} showed that $T(T_X[-1])$ is isomorphic to $(i^*i_*\mathcal{O}_X)^\vee$, where $(-)^\vee$ is the dual.  The map
$$(i^*i_*\mathcal{O}_X)^\vee\rightarrow(i^*j^*j_*i_*\mathcal{O}_X)^\vee=(\Delta^*\Delta_*\mathcal{O}_X)^\vee=\mathcal{RH}om(\Delta_*\mathcal{O}_X,\Delta_*\mathcal{O}_X)$$
is defined by the adjunction $j^*\dashv j_*$.  The composite map
$$\mathrm{Sym}(T_X[-1])\rightarrow T(T_X[-1])\cong(i^*i_*\mathcal{O}_X)^\vee\rightarrow(i^*j^*j_*i_*\mathcal{O}_X)^\vee$$
$$=(\Delta^*\Delta_*\mathcal{O}_X)^\vee=\mathcal{RH}om(\Delta_*\mathcal{O}_X,\Delta_*\mathcal{O}_X) $$

is the sheaf version HKR isomorphism as showed in \cite{AC}.  Taking cohomology on both sides of the equality above, we get the HKR isomorphism
$$I^{HKR}:\mathrm{HT}^*(X)=\bigoplus_{p+q=*}H^p(X,\wedge^q T_X)\rightarrow \mathrm{HH}^*(X).$$

Now it is clear that we have a commutative diagram
$$
\xymatrix{\mathcal{RH}om(\Delta_*\mathcal{O}_X,\Delta_*\mathcal{O}_X)\ar[r] & \mathcal{RH}om(\F,\F)\\
T(T_X[-1])\ar[u]\ar[ur]\\
\mathrm{Sym}(T_X[-1]),\ar[u]\ar[uur]}
$$
which is similar to the Lie algebra diagram in~(\ref{Lie diagram}).
Taking cohomology on the diagram above, we get the diagram
$$
\xymatrix{
\mathrm{HH}^*(X)\ar[r]& \mathrm{Ext}^*(\F,\F)\\
\mathrm{HT}^*(X)=\displaystyle{\bigoplus_{p+q=*}H^p(X,\wedge^qT_X)}\ar[u]^{I^{HKR}}\ar[ur]_{(-)\lrcorner\exp(at_\F)}}
$$
that we start with in Theorem A.
\end{Proof}

\section{The proof of Theorem B}
We use Theorem A to prove Theorem B in this section.

\paragraph  Denote $I^{hkr}(\tilde{\alpha})$ by $\alpha\in\mathrm{HH}^*(X)=\mathrm{Ext}_{X\times X}^*(\mathcal{O}_\Delta,\mathcal{O}_\Delta)$, where $\mathcal{O}_\Delta=\Delta_*\mathcal{O}_X$.  Denote the image of $\alpha$ in $\mathrm{Ext}^*(\F,\F)$ by $\alpha_\F$.  For any vector bundle $\F$ on $X$, C\u{a}ld\u{a}raru and Willerton \cite{CW} defined an abstract Chern character $\mathrm{ch}(\F)$ which lies in the degree zero part of the Hochschild homology $\mathrm{HH}_*(X)=\mathrm{Ext}_{X\times X}^*(S^{-1}_{\Delta},\mathcal{O}_\Delta)$, where $S^{-1}_{\Delta}=\Delta_*(\omega_X^\vee[-\mathrm{dim}X])$.  There is an HKR isomorphism for Hochschild homology
$$I_{HKR}:\mathrm{HH}_*(X)\rightarrow\mathrm{H\Omega}_*(X)=\bigoplus_{q-p=*}H^p(X,\wedge^q\Omega_X).$$
The image of the abstract Chern character under the map $I_{HKR}$ is the usual Chern character of $\F$ \cite{C}.
We need the lemma below.

\begin{Lemma}\label{Lemma}
If $\alpha_\F$ is zero, then $\alpha\circ\mathrm{ch}(\F)$ is zero.  Here $\circ$ is the composition of morphisms in $\mathbf{D}^b(X\times X)$ and $\mathrm{ch}(\F)$ is the abstract Chern character.
\end{Lemma}

\begin{proof}
The proof is known in an email correspondence with Eyal Markman.  Let $\beta$ be any class in $\mathrm{Ext}_{X\times X}^{*}(\mathcal{O}_\Delta,S_\Delta)$, where $S_\Delta=\Delta_*(\omega_X[\mathrm{dim}X])$.  Similar to the definition of the class $\alpha_\F$ associated to $\alpha\in\mathrm{Ext}_{X\times X}^*(\mathcal{O}_\Delta,\mathcal{O}_\Delta)$, we get a class $\beta_\F\in\mathrm{Ext}_X^*(\F,S_X\F)$, where $S_X(-)=\omega_X[\mathrm{dim}X]\otimes-$.  It is shown in \cite{C} that the class $\mathrm{ch}(\F)$ is characterized by the identity
$$\mathrm{Tr}_{X\times X}(\beta\circ\mathrm{ch}(\F))=\mathrm{Tr}_X(\beta_\F).$$
Due to the equality above, we have
$$\mathrm{Tr}_{X\times X}(\gamma\circ\alpha\circ\mathrm{ch}(\F))=\mathrm{Tr}_X((\gamma\circ\alpha)_\F)=\mathrm{Tr}_X(\gamma_\F\circ\alpha_\F)$$
for any $\gamma\in\mathrm{Ext}_{X\times X}^*(\mathcal{O}_\Delta,S_\Delta)$.  The right hand side is zero since we assume that $\alpha_\F$ is zero.  We can conclude that $\alpha\circ\mathrm{ch}(\F)$ is zero because the equality $\mathrm{Tr}_{X\times X}(\gamma\circ\alpha\circ\mathrm{ch}(\F))=0$ holds for any $\gamma$ and $\mathrm{Tr}(-)$ is non-degenerate.
\end{proof}

The two HKR isomorphisms $I_{HKR}$ and $I^{HKR}$ can be twisted by the Todd class.  We denote the resulting twisted isomorphisms  by $I_K$ and $I^K$
$$ I_K: \mathrm{HH}_*(X)\rightarrow \mathrm{H\Omega}_*(X)=\displaystyle{\bigoplus_{q-p=*}H^p(X,\wedge^q\Omega_X)},$$
$$I^K:\mathrm{HT}^*(X)=\bigoplus_{p+q=*}H^p(X,\wedge^qT_X)\rightarrow\mathrm{HH}^*(X).$$

They are given by the formula $I_K=(-\wedge Td(X)^{\frac{1}{2}})\circ I_{HKR}$ and $I^K=I^{HKR}\circ D^{-1}$, where $D^{-1}$ is the inverse of the Duflo operator.

The Mukai vector $v(\F)$ of $\F$ is $I_K(\mathrm{ch}(\F))$ by definition.  There are natural ring structures on $\mathrm{HH}^*(X)$ and $\mathrm{HT}^*(X)$: the product on $\mathrm{HH}^*(X)$ is the Yoneda product, and the product on $\mathrm{HT}^*(X)$ is the wedge product.  Kontsevich \cite{Kont} claimed that the map $I^K$ is a ring isomorphism.  This statement was proved by Calaque and Van den Bergh \cite{CV}.  The Hochschild homology is a module over the Hochschild cohomology and similarly $\mathrm{H\Omega}_*(X)$ is a module over $\mathrm{HT}^*(X)$.  Calaque, Rossi, and Van den Bergh \cite{CRV} proved that the maps $I_K$ and $I^K$ respect the module structures.
\begin{Proof}[Proof of Theorem B.]
The commutative diagram in Theorem A shows that
$$\tilde{\alpha}\lrcorner\exp(at_\F)=\alpha_\F,$$
which is zero under the assumption of Theorem B.  We conclude that $\alpha\circ\mathrm{ch}(\F)$ is zero by Lemma \ref{Lemma}. Since $I_K$ and $I^K$ respect the module structures, we have
$$0=I_K(\alpha\circ\mathrm{ch}(\F))=(I^K)^{-1}(\alpha)\lrcorner\,I_K(\mathrm{ch}(\F))=(I^K)^{-1}(\alpha)\lrcorner \,v(\F).$$  The inverse map of $I^K$ is the composite map
$$\xymatrix{(I^K)^{-1}:\mathrm{HH}^*(X)\ar[r]^{~~~~~~(I^{HKR})^{-1}} &\mathrm{HT}^*(X)\ar[rr]^{D} & & \mathrm{HT}^*(X).}$$
As a consequence $$0=I_{K}(\alpha\circ\mathrm{ch}(\F))=(I^K)^{-1}(\alpha)\lrcorner\,v(\F)=D(\tilde{\alpha})\lrcorner \,v(\F).$$

\end{Proof}

\paragraph{\bf The special case when $\tilde{\alpha}\in H^1(X, T_X)$.}
The result in~(\ref{1.1}) says that $\tilde{\alpha} \lrcorner\,v(\F)$ is zero if $\tilde{\alpha}\lrcorner\exp(at_\F)$ is zero for any $\tilde{\alpha}\in H^1(X, T_X)$.  We end this paper by proving that Theorem B implies the result in~(\ref{1.1}).

From now on let $\tilde{\alpha}$ be an element in $H^1(X, T_X)$.  The only term in $\exp(at_\F)=1+at_\F+\frac{(at_\F)^2}{2!}+\cdots$ that can contract with $\tilde{\alpha}$ is $at_\F$, so $\tilde{\alpha}\lrcorner\exp(at_\F)=\tilde{\alpha}\lrcorner\,at_\F$ in this case.

Choose $\F=O_X$.  We have $\tilde{\alpha}\lrcorner\exp(at_{O_X})=0$.  Therefore $$D(\tilde{\alpha})\lrcorner \,v(O_X)=(Td(X)^{\frac{1}{2}}\lrcorner\,\tilde{\alpha})\lrcorner\,Td(X)^{\frac{1}{2}}=0$$
according to Theorem B.

Expand the Todd class $Td(X)$ as $1+\frac{c_1}{2}+\frac{c_1^2+c_2}{12}+\cdots$, and note that the only term of $(Td(X)^{\frac{1}{2}}\lrcorner\,\tilde{\alpha})\lrcorner\,Td(X)^{\frac{1}{2}}$ in $H^2(X,O_X)$ is $\tilde{\alpha}\lrcorner\,\frac{c_1}{2}$.  Since $(Td(X)^{\frac{1}{2}}\lrcorner\,\tilde{\alpha})\lrcorner \,Td(X)^{\frac{1}{2}}=0$, we can conclude that $\tilde{\alpha}\lrcorner\,c_1$ is zero for any $\tilde{\alpha}\in H^1(X, T_X)$.  The fact that $\tilde{\alpha}\lrcorner\,c_1=0$ for $\tilde{\alpha}\in H^1(X,T_X)$ is also known due to Griffiths.  Consider the first order deformation of $X$ corresponding to $\tilde{\alpha}$.  The vanishing of $\tilde{\alpha}\lrcorner\,c_1$ is equivalent to the class $c_1$ remaining of type $(p,p)$.

The term $\tilde{\alpha}\lrcorner\,\frac{c_1}{4}$ is exactly the difference between $D(\tilde{\alpha}) $ and $\tilde{\alpha} $ because

$$D(\tilde{\alpha})=Td(X)^{\frac{1}{2}} \lrcorner\,\tilde{\alpha}=(1+\frac{c_1}{4}+\cdots)\lrcorner\,\tilde{\alpha}=\tilde{\alpha}+\frac{c_1}{4}\lrcorner\,\tilde{\alpha}+0.  $$

We conclude that $\tilde{\alpha}\lrcorner\,v(\F)$ is zero if and only if $D(\tilde{\alpha})\lrcorner\,v(\F)$ is zero for $\tilde{\alpha}\in H^1(X, T_X)$.

\end{document}